\documentclass[12pt,a4paper]{article}
\usepackage[french]{babel}
\usepackage[T1]{fontenc}
\usepackage[utf8]{inputenc}

\usepackage[backend=bibtex,style=alphabetic,language=french]{biblatex}
\usepackage{csquotes}
\usepackage{amsmath}
\usepackage{amsfonts}
\usepackage{amssymb}
\usepackage{graphicx}
\usepackage{kpfonts}
\usepackage{xcolor}	
\usepackage{xargs}
\usepackage{amsthm}
\usepackage{tikz}
\usepackage[left=1.5cm,right=1.5cm,top=3cm,bottom=3cm]{geometry}
\usepackage{dsfont}

\title{Induction automorphe pour les représentations elliptiques}
\date{}

\makeatletter

\newtheoremstyle{defi}
{12pt}     
{12pt}     
{\normalfont\color{black}}
{}         
{\bf}      
{. }       
{ }   
{\thmname{#1}\thmnumber{\@ifnotempty{#1}{ }\@upn{#2}}%
\thmnote{\ {\the\thm@notefont\bf--- #3.}}}

\newtheoremstyle{prop}
{3pt}     
{3pt}     
{\normalfont\color{black}}
{}         
{\bf}      
{. }       
{ }   
{\thmname{#1}\thmnumber{\@ifnotempty{#1}{ }\@upn{#2}}%
\thmnote{\ {\the\thm@notefont\bf--- #3.}}}

\newtheoremstyle{theoreme}
{3pt}     
{3pt}     
{\normalfont\color{black}}
{}         
{\bf}      
{ }       
{ }   
{\thmname{#1}\thmnumber{\@ifnotempty{#1}{ }\@upn{#2}}%
\thmnote{\ {\the\thm@notefont\bf--- #3.}}}

\newtheoremstyle{lemme}
{3pt}     
{3pt}     
{\normalfont\color{black}}
{}         
{\bf}      
{. }       
{ }   
{\thmname{#1}\thmnumber{\@ifnotempty{#1}{ }\@upn{#2}}%
\thmnote{\ {\the\thm@notefont\bf--- #3.}}}

\newtheoremstyle{preuve}
{3pt}     
{3pt}     
{\normalfont\color{black}}
{}         
{\bf}      
{\ }       
{\newline}   
{\thmname{#1}\thmnumber{\@ifnotempty{#1}{ }\@upn{#2}}%
\thmnote{\ {\the\thm@notefont\bf--- #3.}}}
\makeatother

\begin{document}
\maketitle
\begin{center}
\author{Martin Fatou}
\end{center}

\theoremstyle{defi}
\newtheorem{defi}{Définition}[subsection]
\theoremstyle{prop}
\newtheorem{prop}{Proposition}[subsection]
\theoremstyle{theoreme}
\newtheorem{theo}{Théorème}[subsection]
\theoremstyle{lemme}
\newtheorem{lemme}{Lemme}[section]
\theoremstyle{preuve}
\newtheorem*{dem}{Démonstration}

\section*{Résumé}
Nous étendons l'application de relèvement pour l'induction automorphe définie par une identité de caractères à toutes les représentations elliptiques.

\vskip2mm

\section*{Introduction}
Soit $F$ un corps commutatif localement compact non archimédien, soit $E$ une extension cyclique de $F$ de degré $d$ et soit $m\geq 1$ un entier. D'après le théorème du corps de classes local, l'extension $E$ est définie par un caractère $\kappa: F^{\times}\to \mathbb{C}^{\times}$ tel que $\ker(\kappa)=\text{N}_{E/F}(E^{\times})$, où $\text{N}_{E/F}:E^{\times}\to F^{\times}$ est l'application norme. L'induction automorphe (locale) 
est une application qui associe à une représentation lisse irréductible $\tau$ de $\text{GL}_{m}(E)$ une représentation lisse irréductible $\pi$ de $\text{GL}_{md}(F)$ qui est $\kappa$-stable, 
$i.e.$ isomorphe à $(\kappa\circ \det)\otimes \pi$. Cette application s'exprime par une identité de caractères et correspond, via la correspondance de Langlands locale, à l'induction de $E$ à $F$ des représentations galoisiennes.

L'induction automorphe pour les représentations  génériques unitaires a été démontrée par G. Henniart et R. Herb dans \cite{HH}. 
Nous démontrons ici que cette application existe également pour les représentations elliptiques en utilisant uniquement des arguments locaux.

Pour cela, nous nous inspirons de l'article de A. Badulescu et G. Henniart \cite{BH}, qui concerne le changement de base. Rappelons que le changement de base associe à une représentation lisse irréductible de $\text{GL}_n(F)$ une représentation lisse irréductible $\sigma$-stable de $\text{GL}_n(E)$ où $\sigma$ est un générateur de $\text{Gal}(E/F)$. Tout comme pour l'induction automorphe, l'application de changement de base s'exprime 
par une identité de caractères.

A. Badulescu et G. Henniart démontrent (en particulier) que le changement de base existe pour les représentations elliptiques (Theorem C). Nous suivons de très près leur article.

Nous donnons dans la première section l'identité de caractères définissant l'induction automorphe. Puis nous rappelons les différentes classifications des représentations. L'identité de caractères donnée en section 1 nécessite un opérateur d'entrelacement, c'est pourquoi nous les définissons en section 3. Nous définissons d'abord l'opérateur d'entrelacement d'une induite, puis d'un sous-quotient irréductible et enfin d'un sous-quotient irréductible d'une induite. Nous normalisons ces opérateurs en utilisant les fonctionnelles de Whittaker. En section 4 nous rappelons la construction des représentations elliptiques à partir des représentations essentiellement de carré intégrable. Nous profitons des sections 5 et 6 pour rappeler des résultats déjà établis sur l'induction automorphe: en section 5 l'induction automorphe pour les représentations essentiellement de carré intégrable et en section 6 la compatibilité entre l'induction parabolique et l'induction automorphe. Enfin nous démontrons notre théorème en section 7. Nous montrons que les représentations elliptiques admettent une induite automorphe en exploitant les propriétés des opérateurs d'entrelacement.

\vskip3mm
\textbf{Notations et conventions.}
On note $|.|_F$ et $|.|_E$ les valeurs absolues normalisées de $F$ et $E$.

On note $H$ le groupe $\text{GL}_m(E)$ et $G$ le groupe $\text{GL}_n(F)$ où $n=md$. 

On verra $\kappa$ comme un caractère de $G$, toujours noté $\kappa$, via $\kappa(g)=\kappa(\det g)$ pour $g\in G$. 

Nous ne considérerons que des représentations lisses complexes, i.e. à valeurs dans le groupe des automorphismes d'un espace vectoriel sur $\mathbb{C}$. 
Pour une représentation $\pi$ de $G$, on note $\kappa\pi$ la représentation $(\kappa \circ \det)\otimes \pi$.
\vskip2mm

\section{Définition de l'induction automorphe}

Soit $\tau$ une représentation  irréductible de $H$. 

On définit la notion de $\kappa$-relèvement de $\tau$. 

Pour cela il faut d'abord définir la notion d' intégrales orbitales qui se correspondent puis on définira le $\kappa$-relèvement à l'aide d'une égalité de caractères.

\subsection{Induction parabolique}

Nous ne considérerons dans la suite que des sous-groupes de Levi standards, i.e. des sous-groupes de matrices diagonales par blocs de tailles données. Par exemple pour $G$, si $n_1, \dots, n_k$ sont les tailles des blocs avec $\sum_{i=1}^k n_i=n$, alors $L$, sous-groupe de Levi standard de $G$ associé à $(n_1,\ldots ,n_k)$, est le groupe $\text{GL}_{n_1}(F)\times \text{GL}_{n_2}(F)\times \dots \times \text{GL}_{n_k}(F)$. Nous notons alors $P_L$ le sous-groupe parabolique standard associé, à savoir que $P_L$ est le produit semi-direct $L \rtimes U$ où $U$ est le radical unipotent de $P_L$, c'est-à-dire le groupe des matrices triangulaires supérieures par blocs de tailles $n_1, \dots, n_k$. 

\vskip2mm
Nous noterons alors $\iota_L^G$ l'induction parabolique normalisée de $(L,P_L)$ à $G$.

Si, pour $i=1, \dots, k$,  $\pi_i $ est une représentation de $\text{GL}_{n_i}(F)$, nous notons alors $\pi_1\times \pi_2\times \dots \times \pi_k$ la représentation $\iota_L^G(\pi_1\otimes \pi_2\otimes \dots \otimes \pi_k)$ de $\text{GL}_n(F)$.

\subsection{Facteurs de transfert}
Pour $x\in G$ on écrit $\det\left(T-1+\text{Ad}_G(x)\,\vert\, {\rm Lie}(G)\right)=D_G(x)T^n+\dots$ où $D_G$ est une fonction polynomiale non nulle sur $G$.

On note $G_{\text{reg}}=\lbrace x\in G, D_G(x)\not= 0 \rbrace$ l'ensemble des éléments semisimples réguliers de $G$; c'est encore l'ensemble des éléments de $G$ qui ont $n$ valeurs propres distinctes dans une clôture algébrique de $F$.

On définit de la même manière $D_H$ et $H_{\text{reg}}$. On obtient un plongement de $H$ dans $G$ en fixant une base de $E^m$ en tant que $F$-espace vectoriel. On remarque que $H\cap G_{\text{reg}}\subset H_{\text{reg}}$. 
\vskip3mm

Pour $\gamma, \delta\in H$ soient $c_1,\dots,c_m$ (respectivement $d_1,\dots,d_m$) les valeurs propres de $\gamma$ (respectivement $\delta$) dans une certaine extension de $E$.

On pose:\[r(\gamma, \delta)=\prod_{i,j=1}^m(c_i-d_j).\]
\vskip3mm

Le groupe $\text{Gal}\left(E/F\right)$ agit sur $H$. Soit $\sigma$ un générateur de $\text{Gal}\left(E/F\right)$. Pour $\gamma\in H$ on définit \[\tilde{\Delta}(\gamma)=\prod_{0\leqslant i<j\leqslant d-1}r(\sigma^i\gamma, \sigma^j\gamma).\]

Pour tout $\gamma\in H\cap G_{\text{reg}}$, $\tilde{\Delta}(\gamma)\in E^{\times}$. On sait qu'il existe $e\in E^{\times}$ tel que $e\tilde{\Delta}(\gamma)\in F^{\times}$ pour tout $\gamma\in H\cap G_{\text{reg}}$.
\vskip2mm

Pour $\gamma\in H\cap G_{\text{reg}}$ on pose alors \[\Delta(\gamma)=\kappa\left(e\tilde{\Delta}(\gamma)\right)\](dépend du choix de $e$ et $\sigma$).

On pourra se reporter à \cite{HH} pour les propriétés de ces facteurs de transfert (notamment le paragraphe 4).
\vskip4mm
\subsection{Intégrales orbitales}

Soit $\text{d}g$ une mesure de Haar sur $G$ et $\text{d}h$ sur $H$.

Pour tout $\gamma\in H\cap G_{\text{reg}}$, puisque $\gamma$ est semisimple régulier comme élément de $G$ son centralisateur dans $G$ est un tore $T_{\gamma}$ et ce tore est contenu dans $H$. On fixe sur $T_{\gamma}$ la mesure de Haar $\text{d}t_{\gamma}$ telle que le sous-groupe compact maximal de $T_{\gamma}$ soit de volume 1.

Soient $\frac{\text{d}g}{\text{d}t_{\gamma}}$ et $\frac{\text{d}h}{\text{d}t_{\gamma}}$ les mesures quotient sur $T_{\gamma}\backslash G$ et $T_{\gamma}\backslash H$ respectivement.

On peut maintenant définir les intégrales orbitales.

On note $\mathcal{C}_c^{\infty}\left(G\right)$ l'espace des fonctions complexes sur $G$ qui sont localement constantes et à support compact. Pour $\phi\in \mathcal{C}_c^{\infty}\left(G\right)$ et $\gamma\in G_{\text{reg}}$  on pose\[\Lambda_{\kappa}^G(\phi,\gamma)=\int_{T_{\gamma}\backslash G}\phi(g^{-1}\gamma g)\kappa(g)\frac{\text{d}g}{\text{d}t_{\gamma}}\] si $\gamma$ est tel que $\kappa(g)=1$ pour tout $g\in T_{\gamma}$ (i.e. $T_{\gamma}\subset \ker(\kappa)$), et \[\Lambda_{\kappa}^G(\phi,\gamma)=0\]sinon (observons que si $\gamma\in H\cap G_{\text{reg}}$ on a $\kappa(g)=1$ pour tout $g\in T_{\gamma}$ car $T_\gamma \subset H$ et $\kappa$ est trivial sur $H$).
\vskip3mm

Pour $f\in \mathcal{C}_c^{\infty}\left(H\right)$ et $\gamma\in H_{\text{reg}}$ on pose \[\Lambda^H(f,\gamma)=\int_{T_{\gamma}\backslash H}f(h^{-1}\gamma h)\frac{\text{d}h}{\text{d}t_{\gamma}}.\]

\vskip3mm
On peut alors donner la formulation de l'induction automorphe en termes d'intégrales orbitales.

On dit que $\phi \in \mathcal{C}_c^{\infty}\left(G\right)$ et $f\in \mathcal{C}_c^{\infty}\left(H\right)$ \textit{concordent} ou que $f$ est un \textit{transfert} de $\phi$ si pour tout $\gamma\in H\cap G_{\text{reg}}$, \[\Delta(\gamma)|D_G(\gamma)|_F^{\frac{1}{2}}\Lambda_{\kappa}^G\left(\phi,\gamma\right)=|D_H(\gamma)|_E^{\frac{1}{2}}\Lambda^H\left(f,\gamma\right).\]

\subsection{$\kappa$-relèvement}
Soient $\tau$ une représentation irréductible de $H$, $\pi$ une représentation irréductible de $G$ et $A$ un isomorphisme de $\kappa \pi$ sur $\pi$: $A\circ \kappa \pi(g)=\pi(g)\circ A$ pour tout $g\in G$.
\vskip2mm
Pour $\phi\in \mathcal{C}_c^{\infty}\left(G\right)$, on note $\pi(\phi)$ l'opérateur $v\in V\mapsto \int_G\phi(g)\pi(g)(v)\text{d}g$ où $V$ est l'espace de $\pi$ (de même pour $\tau(f)$).

Puisque $\pi$ et $\tau$ sont admissibles la trace de ces opérateurs est bien définie. 
\vskip2mm

On dit que $\pi$ est un \textit{$\kappa$-rel\`evement} de $\tau$ s'il existe un nombre complexe non nul $c=c(\tau,\pi,A)$ tel que l'on ait \[\text{tr}\left(\pi(\phi)\circ A\right)=c(\tau, \pi, A)\text{tr}\left(\tau(f)\right)\]dès que $\phi\in \mathcal{C}_c^{\infty}\left(G_{\text{reg}}\right)$ et $f\in \mathcal{C}_c^{\infty}\left(H\cap G_{\text{reg}}\right)$ concordent.
\vskip2mm

La notion de $\kappa$-relèvement ne dépend que des classes d'isomorphisme de $\tau$ et $\pi$. 

\subsection{Identité de caractères}
On note $G_0$ le noyau de $\kappa$ vu comme caractère de $G$.

On exploite le fait que la distribution  $\phi\mapsto \text{tr}\left(\pi(\phi)\circ A\right)$ est donnée par une fonction localement constante sur l'ouvert $G_{\text{reg}}$ et de même pour $f\mapsto \text{tr}\left(\tau(f)\right)$. 

Ainsi il existe une fonction $\Theta_{\pi}^A$ localement constante sur $G_{\text{reg}}$ telle que pour tout $\phi\in \mathcal{C}_c^{\infty}\left(G_{\text{reg}}\right)$,\[\text{tr}\left(\pi(\phi)\circ A\right)=\int_{G_{\text{reg}}}\Theta_{\pi}^A(g)\phi(g)\text{d}g\]

et une fonction $\Theta_{\tau}$ localement constante sur $H_{\text{reg}}$ telle que pour tout $f\in \mathcal{C}_c^{\infty}\left(H\cap G_{\text{reg}}\right)$, \[\text{tr}\left(\tau(f)\right)=\int_{H\cap G_{\text{reg}}}\Theta_{\tau}(h)f(h)\text{d}h.\] 

On réécrit alors l'égalité du paragraphe précédent en termes de ces fonctions sur $G$, notamment grâce à la formule d'intégration de Weyl.

\vskip2mm
Soit $\gamma\in H\cap G_{\text{reg}}$, sa classe de conjugaison $\mathcal{O}(\gamma)$ dans $G$ rencontre $H$ en un nombre fini de classes de conjugaison dans $H$. Pour chaque telle classe $C$ dans $H$, on choisit un élément $x_C$ dans $G$ tel que $x_C^{-1}\gamma x_C$ appartienne à $C$, et on note $X(\gamma)$ l'ensemble des éléments $x_C$ pour $C$ parcourant les classes dans $H$ rencontrant $\mathcal{O}(\gamma)$.

\vskip3mm
Alors la condition ``$\pi$ est un $\kappa$-relèvement de $\tau$'' s'écrit avec les égalités suivantes:
\begin{enumerate}
\item pour $\gamma \in H\cap G_{\text{reg}}$, \[|D_G(\gamma)|_F^{\frac{1}{2}}\Theta_{\pi}^A(\gamma)=c(\tau,\pi,A)\sum_{x\in X(\gamma)}\kappa(x^{-1})\Delta(x^{-1}\gamma x)|D_H(x^{-1}\gamma x)|_E^{\frac{1}{2}}\Theta_{\tau}(x^{-1}\gamma x);\]
\item pour $\gamma\in G_{\text{reg}}$ non conjugué à un élément de $H$, \[\Theta_{\pi}^A(\gamma)=0.\]
\end{enumerate}

\section{Classifications}

On énonce les classifications pour $\text{GL}_n(F)$ mais on a les mêmes résultats pour $\text{GL}_m(E)$.

On dispose des classifications suivantes: la classification de Bernstein-Zelevinsky pour les représentations de carré intégrable, la classification de Langlands pour les représentations irréductibles, et la classification de Tadic pour les représentations irréductibles unitaires. 

\subsection{Classification de Bernstein-Zelevinsky}
La classification de Bernstein-Zelevinsky concerne les représentations de carré intégrable.

Soit $\delta$ une représentation irréductible de carré intégrable de $\text{GL}_n(F)$, alors il existe une paire $(k, \rho)$, où $k$ est un diviseur de $n$ et $\rho$ est une représentation irréductible cuspidale unitaire de $\text{GL}_{\frac{n}{k}}(F)$, telle que $\delta$ est isomorphe à l'unique sous-représentation irréductible $Z(\rho, k)$ de $\nu^{\frac{k-1}{2}}\rho\times \nu^{\frac{k-1}{2}-1}\rho\times \dots \times \nu^{-\frac{k-1}{2}}\rho$, où $\nu$ est le caractère de $\text{GL}_n(F)$ égal à la composition de la norme $| \text{ } |_F$ avec l'application déterminant et où l'on induit par rapport au parabolique associé au Levi $\text{GL}_{\frac{n}{k}}(F)\times \dots \times \text{GL}_{\frac{n}{k}}(F)$ ($k$ fois).

L'entier $k$ et la classe d'isomorphisme de $\rho$ sont déterminés par la classe d'isomorphisme de $\delta$.

La représentation $\nu^{\frac{k-1}{2}}\rho\times \nu^{\frac{k-1}{2}-1}\rho\times \dots \times \nu^{-\frac{k-1}{2}}\rho$ a aussi un unique quotient irréductible, son 
quotient de Langlands, que nous définissons au prochain paragraphe.

\vskip3mm
Soit $\delta$ une représentation irréductible essentiellement de carré intégrable de $\text{GL}_n(F)$. Alors il existe un entier $k$ divisant $n$ et une représentation irréductible 
cuspidale $\rho$ de $\text{GL}_{\frac{n}{k}}(F)$ tels que $\delta$ est l'unique sous-représentation irréductible de $\nu^{k-1}\rho\times \nu^{k-2}\rho \times \dots \times \rho$. L'ensemble $\lbrace \rho, \nu\rho, \dots, \nu^{k-1}\rho \rbrace$ s'appelle le \textit{segment de Zelevinsky} de $\delta$, l'entier $k$ est sa longueur. 

Notons que $\delta$ est de carré intégrable (i.e. unitaire) si et seulement si $\rho'=\nu^{\frac{k-1}{2}}\rho$ est unitaire, auquel cas l'unique sous-représentation irréductible de $\nu^{k-1}\rho\times \nu^{k-2}\rho \times \dots \times \rho$ est $Z(\rho',k)$.

\subsection{Classification de Langlands}
La classification de Langlands exprime une représentation irréductible en fonction de représentations tempérées.

Soit $n\geqslant 1$ un entier et soit $n=\sum_{i=1}^kn_i$ une partition de $n$ pour des entiers $n_i\geqslant 1$.

Soient $\alpha_1,\dots ,\alpha_k$ des nombres réels tels que $\alpha_1>\alpha_2>\dots >\alpha_k$.

Soient $\tau_1,\dots,\tau_k$ des représentations irréductibles tempérées des groupes $\text{GL}_{n_i}(F)$.

Alors la représentation $\nu^{\alpha_1}\tau_1\times \nu^{\alpha_2}\tau_2\times \dots \times \nu^{\alpha_k}\tau_k$ a un unique quotient irréductible, appelé le \textit{quotient de Langlands} et noté $L(\nu^{\alpha_1}\tau_1, \nu^{\alpha_2}\tau_2, \dots , \nu^{\alpha_k}\tau_k)$.

La classification de Langlands énonce alors que toute représentation irréductible $\pi$ de $\text{GL}_n(F)$ est isomorphe à un tel $L(\nu^{\alpha_1}\tau_1, \nu^{\alpha_2}\tau_2, \dots, \nu^{\alpha_k}\tau_k)$ où $k$, les réels $\alpha_1,\dots ,\alpha_k$ et les classes d'isomorphisme des représentations irréductibles tempérées $\tau_1, \dots, \tau_k$ sont déterminés par la classe d'isomorphisme de $\pi$.

\subsection{Représentations de Speh}
Soit $\tau$ une représentation irréductible tempérée de ${\rm GL}_n(F)$ et $k\geqslant 1$ un entier.

On note alors $u(\tau, k)$ la représentation $L(\nu^{\frac{k-1}{2}}\tau, \nu^{\frac{k-1}{2}-1}\tau,\dots ,\nu^{-\frac{k-1}{2}}\tau)$.

Lorsque $\tau$ est de carré intégrable, $u(\tau,k)$ est appelée \textit{représentation de Speh}.

Si $\alpha\in ]0,\frac{1}{2}[$, on note $\pi(u(\tau,k),\alpha)$ la représentation $\nu^{\alpha}u(\tau,k)\times \nu^{-\alpha}u(\tau,k)$ qui est irréductible.

\subsection{Classification de Tadic}

Soit $\mathcal{U}$ l'ensemble des classes d'isomorphisme de toutes les représentations de la forme $u(\tau,k)$ et $\pi(u(\tau,k),\alpha)$ où $k\geqslant 1$ est un entier, $\tau$ est une représentation de carré intégrable de $\text{GL}_r(F)$ et $\alpha\in ]0,\frac{1}{2}[$.

Alors tout produit d'éléments de $\mathcal{U}$ est irréductible et unitaire. Inversement, toute représentation irréductible unitaire de $\text{GL}_n(F)$ est un produit d'éléments de $\mathcal{U}$ et les facteurs du produit sont déterminés par cette représentation.

\section{Opérateurs d'entrelacement}

Pour $\pi$ une représentation d'un groupe $G$ on notera $V_{\pi}$ l'espace vectoriel associé.

Dans cette partie, nous définissons d'abord l'opérateur d'entrelacement d'une induite $\kappa$-stable à partir de l'opérateur d'entrelacement de la représentation induisante elle aussi 
supposée $\kappa$-stable. Puis nous définissons l'opérateur d'entrelacement d'un 
sous-quotient irréductible de cette induite grâce à la propriété de multiplicité $1$. Ensuite nous mélangeons ces deux propriétés pour obtenir la propriété d'induction parabolique et de multiplicité 1. Enfin nous normalisons ces opérateurs d'entrelacement.

\subsection{Définition de l'opérateur d'entrelacement d'une induite}

Soit $\tau$ une représentation $\kappa$-stable d'un Levi (standard) $L$ de $G$ et soit $P=P_L$. On note $B: \kappa\tau \to \tau $ un opérateur d'entrelacement, i.e. un $L$-isomorphisme entre $\kappa\tau= (\kappa\circ \det)\otimes \tau$ et $\tau$.

\textbf{Vocabulaire.} Pour une représentation $\pi_0$, nous appellerons \textit{$\kappa$-opérateur sur $\pi_0$} un opérateur d'entrelacement entre $\kappa\pi_0$ et $\pi_0$.  
\vskip1mm

Alors montrons que $\pi=\iota_{L}^G(\tau)$ est $\kappa$-stable et que  $A: f\in V_{\pi}\mapsto \left(g\mapsto \kappa(g)B(f(g))\right)$ est un isomorphisme de $\kappa\pi$ sur $\pi$.
\vskip1mm
La représentation $\kappa\pi$ agit sur le même espace $V_{\pi}$ que $\pi$ et l'action est donnée par $\kappa\pi(g)=\kappa\circ\det(g)\pi(g)$ pour $g\in G$, $\pi$ agissant par translations à droite sur $V_{\pi}$: pour $f\in V_{\pi}$, $g,g'\in G$, $\pi(g)(f)(g')=f(gg')$.
\begin{itemize}

\item Vérifions que pour $f\in V_{\pi}$ on a bien $Af \in V_{\pi}$.

Soient $p\in P$ et $g\in G$. Comme $f\in V_{\pi}$ on sait que $f(pg)=\delta^{1/2}(p)\tau(p)(f(g))$, d'où \[\left(Af\right)(pg)=\kappa(pg)B\left(f(pg)\right)=\kappa(p)\kappa(g)B\left(\delta^{1/2}(p)\tau(p)(f(g))\right)=\delta^{1/2}(p)\kappa(g)B\left(\kappa\tau(p)(f(g))\right). \]
Or, $B$ est un opérateur d'entrelacement entre $\kappa\tau$ et $\tau$, on obtient donc  \[\left(Af\right)(pg)=\delta^{1/2}(p)\kappa(g)\tau(p)B\left(f(g)\right)=\delta^{1/2}(p)\tau(p)\left(Af\right)(g)\]
i.e. $Af\in V_{\pi}$.

\item Vérifions maintenant que $A$ est bien un opérateur d'entrelacement entre $\kappa\pi$ et $\pi$, i.e. pour tout $g\in G$, $A\circ \kappa\pi(g)=\pi(g)\circ A$.

Soient donc $g\in G, f\in V_{\pi}$ et $g'\in G$. On a 

\begin{align*}
\left(A\circ\kappa\pi(g)\right)(f)(g')=A\left(\kappa\pi(g)(f)\right)(g')=\kappa(g')B\left(\kappa\pi(g)(f)(g')\right)&=\kappa(g')B\left(\kappa(g)f(gg')\right) \\
&=\kappa(gg')B\left(f(gg')\right) \\
&=\left(Af\right)(gg') \\
&=\left(\pi(g)\circ A\right)(f)(g').
\end{align*}

Donc $A$ est bien un opérateur d'entrelacement entre $\kappa\pi$ et $\pi$.

\end{itemize}

\subsection{Propriété de multiplicité 1}
On reprend en l'adaptant le paragraphe 2.2 de l'article de Badulescu-Henniart \cite{BH}. On considère ici un groupe localement profini $G$, 
un caractère $\kappa$ (i.e. un homomorphisme continu dans ${\mathbb C}^\times$) de $G$ et une représentation (complexe, lisse) $\pi$ de $G$. 
\begin{itemize}
\item Pour le changement de base on prend un isomorphisme de $\pi$ sur $\pi^{\sigma}$ alors qu'ici on prend un isomorphisme de $\kappa\pi$ sur $\pi$.
\item Tout ce qui est dans l'Appendix de \cite{BH} peut être repris: la première partie ne concerne que des résultats d'algèbre générale, 
la deuxième partie ("Group with automorphism") est à adapter avec $\kappa$ au lieu de $\sigma$ compte-tenu des propriétés: 
\begin{itemize}
\item un sous-espace de $V_{\pi}$ est stable par $\pi$ si, et seulement si, il est stable par $\kappa\pi$ (rappelons que $\kappa\pi$ opère naturellement sur $V_\pi$);
\item si $U,W$ sont des sous-espaces stables de $V_{\pi}$ tels que $W\subset U$, on note $\pi_U$ la sous-représentation de $\pi$ dans $U$, et $\pi_{U/W}$ la représentation quotient de $\pi_U$ dans $U/W$ 
induite par $\pi$. On a alors $\kappa\left(\pi_U\right)=\left(\kappa\pi\right)_U$ et $\kappa\left(\pi_{U/W}\right)=\left(\kappa\pi\right)_{U/W}$.
\end{itemize}
\end{itemize}

\vskip1mm

On peut donc appliquer la propriété de \cite[Appendix]{BH}, appelée \textit{"propriété de multiplicité 1"} que l'on rappelle ci-dessous.

Cette propriété concerne le lien entre les isomorphismes d'une représentation et les sous-quotients irréductibles de cette représentation.

On suppose que $\pi$ est de longueur finie et $\kappa$-stable. On fixe $f: \kappa\pi\to \pi$ un $G$-isomorphisme. Soit $\pi_0$ un sous-quotient irréductible de $\pi$, supposé $\kappa$-stable.

On suppose de plus que $\pi_0$ est de multiplicité $1$ dans $\pi$. Il existe une paire $(U,W)$ de sous-espaces stables de $V_\pi$ avec $W\subset U$ telle que $\pi_0 \simeq \pi_{U/W}$. Si de plus $U$ est maximal pour cette propriété, 
ce que l'on suppose, alors la paire $(U,W)$ est déterminée de manière unique [BH, proposition 7.1, (b)]. 
On fixe un $G$-isomorphisme $\phi: \pi_0 \simeq \pi_{U/W}$.

L'application $f$ induit par passage au quotient un $G$-isomorphisme $\overline{f}:\kappa\pi_{U/W}\to \pi_{U/W}$. Alors on obtient un opérateur d'entrelacement \[\phi^{-1}\overline{f}\phi: \kappa\pi_0\to \pi_0\]qui ne dépend pas du choix de $\phi$ (lemme de Schur).

On dit que l'opérateur $\phi^{-1}\overline{f}\phi$ est le $\kappa$-opérateur sur $\pi_0$ obtenu à partir de $f$ par la propriété de multiplicité 1.

\subsection{Induction parabolique et multiplicité 1}

Reprenons les hypothèses et les notations de 3.1. On peut noter $B_{\kappa}(\pi)$ l'opérateur d'entrelacement $A$ défini en loc.~cit., qui est l'équivalent de l'opérateur $I_s(\pi)$ du paragraphe 2.2 de \cite{BH}.

Pour rappel, pour $f\in V_{\pi}$, $B_{\kappa}(\pi)(f)$ est donné par $B_{\kappa}(\pi)(f)(g)=\kappa(g)B\left(f(g)\right)$ pour $g\in G$.

Nous énonçons ici la \textit{propriété d'induction parabolique et de multiplicité 1} qui consiste à mixer les deux constructions précédentes 
et donc à construire un opérateur d'entrelacement sur un sous-quotient irréductible de multiplicité $1$ d'une représentation induite parabolique.

Plus précisément, soit $\tau$ une représentation $\kappa$-stable d'un Levi $L$, soit $B:\kappa\tau\to\tau$ un opérateur d'entrelacement et soit $\pi=\iota_{L}^{G}(\tau)$. Nous rencontrerons souvent la situation où $\pi$ a un sous-quotient $\pi_0$ irréductible de multiplicité $1$ et $\kappa$-stable. Alors, d'après 3.2 le $\kappa$-opérateur $B_{\kappa}(\pi)$ sur $\pi$ obtenu à partir de $B$ par induction parabolique (cf. ci-dessus) induit par la propriété de multiplicité 1 un opérateur $B_{\kappa}(\pi_0)$ sur $\pi_0$ qui est bien défini, i.e. ne dépend pas de la manière dont on réalise $\pi_0$ comme sous-quotient de $\pi$.

\textbf{Définition}  On dit que $B_{\kappa}(\pi_0)$ est le $\kappa$-opérateur sur $\pi_0$ obtenu à partir de $B$ par la propriété d'induction parabolique et de multiplicité 1.

\subparagraph{Comportement des $\kappa$-opérateurs avec l'induction parabolique}
\vskip1mm

On démontre que la proposition 2.1 de \cite{BH} est toujours valable pour les $\kappa$-opérateurs.

\vskip3mm

Soit $L'$ un sous-groupe de Levi de $G$ tel que $L\subset L'$. 

On pose \[V'=\lbrace f:L'\to W \text{ lisse }, f(pg)=\delta^{1/2}(p)\tau(p)f(g) \forall g\in L', p\in P_L\cap L'\rbrace\] et $\tau'=\iota_{L}^{L'}\tau$ la représentation par translations à droite de $L'$ dans $V'$.

\vskip2mm

On induit encore, la représentation $\iota_{L'}^G\tau'$ est la représentation par translations à droite de $G$ dans $V''$ où \[V''=\lbrace f:G\to V' \text{ lisse }, f(pg)=\delta^{1/2}(p)\tau'(p)f(g) \forall g\in G, p\in P_{L'}\rbrace.\]

\vskip2mm
On sait alors (transitivité du foncteur induction parabolique) qu'il existe un isomorphisme $h:V''\to V_{\pi}$ entre $\iota_{L'}^G\tau'$ et $\pi$ défini pour $f\in V''$ par $h(f)=\left(g\in G \mapsto f(g)(1)\right)$.

\vskip3mm
On a alors la proposition suivante qui nous dit principalement (deuxième point) que le $\kappa$-opérateur sur $\pi_0$ défini plus haut ne "dépend pas" de la réalisation de l'induite parabolique dont $\pi_0$ en est un sous-quotient.

\prop \label{transitivité} 
On suppose que $\tau$ est $\kappa$-stable et que $B$ entrelace $\kappa\tau$ et $\tau$.
\begin{enumerate}
\item On a $h\circ \left(B_{\kappa}(\tau')\right)_{\kappa}(\pi)=B_{\kappa}(\pi)\circ h$.
\item Soit $\pi_0$ un sous-quotient $\kappa$-stable irréductible de $\pi$ de multiplicité 1. Soit $\tau_0'$ le sous-quotient irréductible de $\tau'$ tel que $\pi_0$ soit un sous-quotient de $\iota_{L'}^{G_E}(\tau_0')$. Si $\tau_0'$ est $\kappa$-stable on a: \[B_{\kappa}(\pi_0)=\left(B_{\kappa}(\tau_0')\right)_{\kappa}(\pi_0).\]
\end{enumerate}

\textit{Démonstration}
\begin{enumerate}

\item Soient donc $f\in V''$ et $g\in G$. On a \[\left(B_{\kappa}(\pi)\circ h\right)(f)(g)=B_{\kappa}(\pi)\left(h(f)\right)(g)=\kappa(g)B\left(h(f)(g)\right)=\kappa(g)B\left(f(g)(1)\right).\]

D'autre part,\[\left(h\circ \left(B_{\kappa}(\tau')\right)_{\kappa}(\pi)\right)(f)(g)=\left(\left(B_{\kappa}(\tau')\right)_{\kappa}(\pi)(f)\right)(g)(1)=\kappa(g)B_{\kappa}(\tau')\left(f(g)\right)(1)=\kappa(g)B\left(f(g)(1)\right).\]

\item Soit \[0\to W\to U \to \tau_0'\to 0\]une suite exacte de représentations, où $(U,W)$ est la paire maximale de sous-représentations de $\tau'$ telle que $\tau_0'\simeq U/W$. D'après \cite[prop.~7.1\, (c)]{BH}, $U$ et $W$ sont stables par $B_{\kappa}(\tau')$. 

Le foncteur induction parabolique $\iota_{L'}^G$ est exact et on obtient la suite exacte de $G$-modules:
\[0\to \iota_{L'}^{G}W\to \iota_{L'}^{G}U \overset{F}{\to}\iota_{L'}^{G}\tau_0'\to 0.\]

$\pi_0$ est un sous-quotient de $\iota_{L'}^{G}\tau_0'$ de multiplicité 1, soit $(u,w)$ la paire maximale de sous-représentations de $\iota_{L'}^{G}\tau_0'$ telle que $\pi_0 \simeq u/w$.

On a donc une chaîne d'inclusions \[\iota_{L'}^{G}W\subset F^{-1}(w)\subset F^{-1}(u)\subset \iota_{L'}^{G}U\] telle que l'isomorphisme (déduit de $F$)
\[\iota_{L'}^{G}U/\iota_{L'}^{G}W\simeq \iota_{L'}^{G}(\tau_0')\]envoie $F^{-1}(u)/\iota_{L'}^{G}W$ sur $u$ et $F^{-1}(w)/\iota_{L'}^{G}W$ sur $w$.

Par le point 1, en notant $\tilde{B}$ l'opérateur obtenu sur $U$ par restriction de $B_{\kappa}(\tau')$, on sait que l'opérateur sur $\iota_{L'}^{G}U$ obtenu par restriction de $B_{\kappa}(\pi)$ est égal à $\tilde{B}_{\kappa}\left(\iota_{L'}^{G}U\right)$.

Or, $F^{-1}(u)$ est le sous-module maximal de $\iota_{L'}^{G}U$ admettant $\pi_0$ comme quotient.

Donc, les deux opérateurs, construits grâce à la propriété de multiplicité 1 en utilisant les deux façons de voir $\pi_0$ comme un quotient, coïncident, i.e. \[B_{\kappa}(\pi_0)=\left(B_{\kappa}(\tau_0')\right)_{\kappa}(\pi_0).\]

\end{enumerate}

\subsection{Normalisation}

Enfin il reste à traiter l'équivalent de l'"opérateur de $\sigma$-entrelacement normalisé".

On part d'une représentation $\kappa$-stable et on veut normaliser l'opérateur d'entrelacement $A$. Pour cela on traite d'abord le cas des représentations génériques puis on obtient le cas général grâce à la classification de Langlands.

Rappelons ce qu'est une représentation générique. 

On fixe un caractère additif non trivial  $\psi$ de $F$. On obtient un caractère $\theta=\theta_{\psi}$ du sous-groupe unipotent supérieur $U$ de $G$ via:\[\theta (u)=\psi\left(\sum_{i=1}^{n-1}u_{i,i+1}\right) \text{ pour } u=(u_{i,j})\in U.\]

Soit $\pi$ une représentation irréductible de $G$. On dit que $\pi$ est \textit{générique} s'il existe une forme linéaire non nulle $\lambda$ sur l'espace $V_{\pi}$ de $\pi$ telle que l'on ait $\lambda\left(\pi(u)(v)\right)=\theta(u)\lambda(v)$ pour $u\in U$ et $v\in V_{\pi}$.

Cette existence ne dépend pas du choix de $\psi$, et $\lambda$ est unique à un scalaire près, on l'appelle \textit{fonctionnelle de Whittaker} pour $\pi$ relative à $\psi$. On note $\mathcal{W}(\pi,\psi)$ leur ensemble.

\vskip2mm
Soit $\lambda$ une fonctionnelle de Whittaker pour $\pi$ relative à $\psi$. Alors pour $u\in U$ et $v\in V_{\pi}$, \[\lambda\left(\kappa\pi(u)(v)\right)=\kappa\circ\det(u)\lambda\left(\pi(u)(v)\right)=\theta(u)\lambda(v)\]car $\det(u)=1$.

Donc $\lambda$ est également une fonctionnelle de Whittaker pour $\kappa\pi$ relative à $\psi$.

Si $\pi$ est $\kappa$-stable, en notant $A$ l'isomorphisme entre $\kappa\pi$ et $\pi$, on normalise $A$ en imposant $\lambda\circ A=\lambda$ pour toute fonctionnelle de Whittaker $\lambda$. On note $A^{\text{gén}}(\pi,\psi)$ cet opérateur d'entrelacement normalisé, on a donc \[\lambda\circ A^{\text{gén}}(\pi, \psi)=\lambda.\]

Contrairement au cas du changement de base, cet opérateur $A^{\text{gén}}(\pi, \psi)$ dépend du choix de $\psi$.

Si $a\in F^{\times}$ et si on note $\psi^a$ le caractère $x\in F\mapsto \psi(ax)$ alors on a \[A^{\text{gén}}(\pi,\psi^a)=\kappa(t_a)^{-1}A^{\text{gén}}(\pi,\psi)\]où $t_a=\text{diag}(a^{n-1},a^{n-2},\dots,a,1)$.

En effet, on a un isomorphisme \[\lambda\in \mathcal{W}(\pi, \psi)\to \lambda\circ\pi(t_a)\in \mathcal{W}(\pi, \psi^a).\]Par définition de $A^{\text{gén}}(\pi,\psi^a)$ on a l'égalité \[\lambda'\circ A^{\text{gén}}(\pi,\psi^a)=\lambda' \text{ pour tout } \lambda'\in \mathcal{W}(\pi,\psi^a)\]donc en particulier, grâce à l'isomorphisme ci-dessus, pour tout $\lambda\in \mathcal{W}(\pi,\psi)$ on a \[\lambda\circ\pi(t_a)\circ A^{\text{gén}}(\pi,\psi^a)=\lambda\circ\pi(t_a).\]

D'où \[\pi(t_a)\circ A^{\text{gén}}(\pi,\psi^a)\circ \pi(t_a)^{-1}=A^{\text{gén}}(\pi,\psi)\] et donc, comme $A^{\text{gén}}(\pi,\psi)$ entrelace $\kappa\pi$ et $\pi$, 
\begin{align*}
A^{\text{gén}}(\pi,\psi^a)&=\pi(t_a)^{-1}\circ A^{\text{gén}}(\pi,\psi)\circ \pi(t_a) \\
&=\pi(t_a)^{-1}\kappa(t_a)^{-1}\circ A^{\text{gén}}(\pi,\psi)\circ \kappa\pi(t_a) \\
&=\kappa(t_a)^{-1}\pi(t_a)^{-1}\pi(t_a) A^{\text{gén}}(\pi,\psi)\\
&=\kappa(t_a)^{-1}A^{\text{gén}}(\pi,\psi).
\end{align*}

On note $A^{\text{gén}}_{\pi}$ l'opérateur $A^{\text{gén}}(\pi,\psi)$ lorsque le caractère $\psi$ est sous-entendu ou si son choix n'est pas important pour les résultats en question.

\subsubsection*{Comportement de $A^{\text{gén}}_{\pi}$ avec les isomorphismes}
\vskip1mm 

Soit $\pi$ une représentation irréductible générique de $G$. Si $\pi'$ est une représentation isomorphe à $\pi$ et si on note $\phi: \pi\to \pi'$ un isomorphisme alors $\lambda'\mapsto \lambda'\circ \phi$ est un isomorphisme de $\mathcal{W}(\pi', \psi)$ sur $\mathcal{W}(\pi, \psi)$. 

Donc $A^{\text{gén}}_{\pi'}\circ \phi=\phi\circ A^{\text{gén}}_{\pi}$.

\subsubsection*{Définition de $A(\pi, \psi)$ pour $\pi$ irréductible et $\kappa$-stable} 
\vskip1mm

On s'appuie sur la classification de Langlands.

Soit $\pi$ une représentation irréductible $\kappa$-stable de $G$. On sait que $\pi$ est isomorphe à un quotient de Langlands $L(\pi_1,\dots,\pi_r)$ où les $\pi_i$ sont essentiellement tempérées (et donc génériques). Par unicité du quotient de Langlands, les $\pi_i$ sont également $\kappa$-stables. On a donc un opérateur d'entrelacement normalisé $A^{\text{gén}}_{\pi_i}$ entre $\kappa\pi_i$ et $\pi_i$, et on obtient un opérateur d'entrelacement $B=A^{\text{gén}}_{\pi_1}\otimes\dots \otimes A^{\text{gén}}_{\pi_r}$ entre les représentations $\kappa\pi_1\otimes \dots \otimes \kappa\pi_r$ et $\pi_1\otimes \dots \otimes \pi_r$ de $L$, où $L=\text{GL}(n_1,F)\times \dots \times \text{GL}(n_r,F)$ est le sous-groupe de Levi de $G$ sur lequel vit la représentation $\pi_1\otimes \dots \otimes \pi_r$.

Cet opérateur donne par induction parabolique un opérateur d'entrelacement $A=\iota_{L}^{G}(B)$ entre $\kappa\Sigma$ et $\Sigma$ où $\Sigma=\pi_1\times \dots \times \pi_r$.

On note $\overline{\Sigma}=L(\pi_1,\dots,\pi_r)$ l'unique quotient irréductible de $\Sigma$. Puisque $\Sigma$ est $\kappa$-stable, $\kappa\overline{\Sigma}\simeq\overline{\Sigma}$ et $\kappa\overline{\Sigma}$ est l'unique quotient irréductible de $\kappa\Sigma$.

Par passage au quotient, l'opérateur $A$ induit un isomorphisme $\overline{A}$ entre $\kappa\overline{\Sigma}$ et $\overline{\Sigma}$.

Par construction $\pi\simeq\overline{\Sigma}$, il existe donc un morphisme surjectif de $G$-modules $f:V_{\Sigma}\to V_{\pi}$, qui se factorise en un isomorphisme \[\overline{f}\in \text{Isom}_{G}(\overline{\Sigma},\pi)=\text{Isom}_{G}(\kappa\overline{\Sigma},\kappa\pi)\]qui ne dépend pas de $f$ à multiplication près par une constante (lemme de Schur).

\vskip2mm
On définit alors l'opérateur $A(\pi,\psi)$ par \[A(\pi,\psi):=\overline{f}\circ \overline{A}\circ \overline{f}^{-1},\]
opérateur qui ne dépend pas du choix de $f$.

\subsubsection*{$A_{\pi}$ est compatible avec les isomorphismes}
\vskip1mm

Si $\pi'$ est une autre représentation lisse irréductible $\kappa$-stable de $G$ telle que $\pi'\simeq \pi$ et si $\phi\in \text{Isom}_{G}(\pi, \pi')$ alors $\overline{f'}=\phi\circ \overline{f}:V_{\overline{\Sigma}}\to V_{\pi'}$ est un isomorphisme surjectif de $G$-modules et, d'après le point précédent \[A(\pi',\psi)=\overline{f'}\circ \overline{A}\circ (\overline{f'})^{-1}=\phi\circ A(\pi,\psi)\circ \phi^{-1}.\]
D'où \[A(\pi',\psi)\circ \phi=\phi\circ A(\pi, \psi).\]

\subsubsection*{$A_{\pi}$ bien défini} 
\vskip1mm

Montrons que la définition ci-dessus est bien correcte dans le sens où l'opérateur $A(\pi, \psi)$ coïncide avec $A^{\text{gén}}(\pi, \psi)$ lorsque $\pi$ est générique.

\vskip2mm
Soit donc $\pi$ générique. On écrit $\pi$ comme un quotient de Langlands $L(\pi_1,\dots, \pi_r)$. Comme $\pi$ est générique, le produit $\pi_1\times \dots \times \pi_r$ est irréductible et $L(\pi_1,\dots, \pi_r)=\pi_1\times \dots \times \pi_r$. Par compatibilité avec les isomorphismes (point précédent) on peut en fait supposer que $\pi=\pi_1\times \dots \times \pi_r$.

\vskip2mm
Il s'agit de vérifier que $A(\pi, \psi)$ vérifie \[\Lambda\circ A(\pi,\psi)=\Lambda\]pour $\Lambda$ une fonctionnelle de Whittaker pour $\pi$ relative à $\psi$. On commence par construire une telle fonctionnelle de Whittaker.

Soit, pour chaque $i\in \lbrace 1, \dots, r\rbrace$,  $\lambda_i\in \mathcal{W}(\pi_i,\psi)$. D'après \cite[formula (2) chapter 3]{JS} on a alors une fonctionnelle de Whittaker $\Lambda$ sur $\pi_1\times \dots \times \pi_r$ donnée par \[\Lambda(f)=\int_{U}\lambda(f(u))\overline{\Theta_{\psi}(u)}\text{d}u\]où $\lambda=\lambda_1\otimes\lambda_2\otimes \dots \otimes \lambda_k$, $f$ est une fonction dans l'espace de $\pi_1\times \dots \times \pi_k$ et $\text{d}u$ est une mesure de Haar sur $U$; l'intégrale étant toujours convergente d'après \cite{JS}.

\vskip2mm
Alors, pour $u\in U$ on a \[\lambda\left(A(\pi,\psi)(f)(u)\right)=\lambda\left(\kappa(u)A_{\pi_1}^{\text{gén}}\otimes\dots\otimes A_{\pi_r}^{\text{gén}}\left(f(u)\right)\right).\]

Or, pour $u\in U$, $\kappa(u)=\kappa\circ\det(u)=1$ et \[\lambda\circ A_{\pi_1}^{\text{gén}}\otimes\dots\otimes A_{\pi_r}^{\text{gén}}=\lambda\]par définition des $A_{\pi_i}^{\text{gén}}$ et de $\lambda=\lambda_1\otimes \dots \otimes \lambda_r$.

D'où \[\lambda\left(A(\pi,\psi)(f)(u)\right)=\lambda(f(u))\]et donc \[\Lambda\circ A(\pi, \psi)=\Lambda,\]ce que l'on voulait.

\subsubsection*{Compatibilité entre l'induction parabolique et les $\kappa$-opérateurs normalisés} 

\vskip1mm

La proposition suivante exprime la compatibilité entre l'induction parabolique et l'opérateur de $\kappa$-entrelacement normalisé.

\prop \label{compatibilitéinductionnormalisation} Soit $L$ un sous-groupe de Levi standard de $G$, $\tau$ une représentation générique $\kappa$-stable de $L$ et $A_{\tau}=A^{\text{gén}}(\tau, \psi)$ l'opérateur de $\kappa$-entrelacement normalisé de $\tau$.
Alors $\iota_L^{G}(\tau)$ a un unique sous-quotient irréductible générique $\pi_0$ qui est $\kappa$-stable. Si on note $A_{\tau,\kappa}(\pi_0)$ l'opérateur sur $\pi_0$ obtenu à partir de $A_\tau$ 
par la propriété d'induction parabolique et de multiplicité $1$ (3.3.1), alors $A_{\tau,\kappa}(\pi_0)=A^{\text{gén}}(\pi_0, \psi)$.

\vskip3mm

\dem

Pour $\lambda$ une fonctionnelle de Whittaker pour $\tau$ on associe une fonctionnelle de Whittaker pour $\iota_L^{G}(\tau)$ via $\lambda\mapsto (f\in V_{\iota_L^{G}(\tau)}\mapsto \lambda \circ f)$. On sait que $\tau$ est générique, donc  $\iota_L^{G}(\tau)$ a une unique droite de fonctionnelles de Whittaker d'après \cite{JS}  et donc il y a un unique sous-quotient irréductible $\pi_0$ avec des fonctionnelles de Whittaker non nulles, i.e. $\pi_0$ générique et donc $\kappa\pi_0$ générique. 

On sait que $\pi= \iota_L^{G}(\tau)$ est $\kappa$-stable donc par la propriété de multiplicité 1 on obtient que $\pi_0$ est $\kappa$-stable.

On note $\pi_0=U/V$ avec $V\subset U \subset V_\pi$ et $U$ maximal. Alors $U$ et $V$ sont stables par $A_{\tau,\kappa}(\pi)$ qui induit donc par passage aux quotients un opérateur $A_{\tau, \kappa}(\pi_0)$ sur $\pi_0$.

Si $\Lambda$ est une fonctionnelle de Whittaker non nulle pour $\pi$ alors elle induit par restriction une fonctionnelle $\Lambda_U$ sur $U$. De la même manière que dans la preuve du point précédent, 
l'opérateur $A _{\tau,\kappa}(\pi)$ fixe $\Lambda$ et donc sa restriction à $U$ fixe $\Lambda_U$. Donc $A_{\tau, \kappa}(\pi_0)=A^\text{gén}(\pi_0)$. \hfill \qed

\section{Construction de représentations elliptiques}

Rappelons la construction de représentations elliptiques. Notre théorème traitant des représentations elliptiques de $H$, nous donnons ici la construction de représentations elliptiques de $H$ pour conserver les mêmes notations dans le paragraphe 7. Partons d'une représentation essentiellement de carré intégrable $\tau_E$ de $H$ à laquelle on associe des représentations elliptiques comme suit.

D'après la classification de Bernstein-Zelevinsky (voir 2.1), il existe un entier $k$ divisant $m$ et une représentation cuspidale $\rho_E$ de $\text{GL}_{\frac{m}{k}}(E)$ tels que $\tau_E$ se réalise comme l'unique sous-représentation irréductible de l'induite parabolique \[\nu^{k-1}\rho_E\times\nu^{k-2}\rho_E\times \dots \times \rho_E. \]

Pour $I$ un sous-ensemble de $\mathcal{K}=\lbrace 1,\dots, k-1\rbrace$, on définit  un sous-groupe de Levi $L_{E,I}$ de $H$ contenant $L_E=\text{GL}_{\frac{m}{k}}(E)\times \dots \times \text{GL}_{\frac{m}{k}}(E)$ 
de la manière suivante: si $I$ est le complémentaire de $\lbrace n_1,n_1+n_2, \dots, n_1+n_2+\dots +n_{t-1} \rbrace$ dans $\lbrace 1, \dots, k-1 \rbrace$, alors on pose  $$L_{E,I}=\text{GL}_{n_1\frac{m}{k}}(E)\times \dots \times \text{GL}_{n_t\frac{m}{k}}(E)$$
où $n_t$ est tel que $n_1+n_2+\dots +n_t=k$. On a alors $$L_{E,I}\subset L_{E,J} \text{ si } I\subset J$$ et en particulier $L_{E,\emptyset}=L_E$ et $L_{E, \mathcal{K}}=H$.

Pour chaque sous-ensemble $I$ de $\mathcal{K}$ on note:
\begin{itemize}
\item $\tau_{E,I}$ l'unique sous-représentation irréductible de $\iota_{L_E}^{L_{E,I}}(\nu_E^{k-1}\rho_E \otimes \dots \otimes \rho_E)$;
\item $\pi_{E,I}$ le quotient de Langlands, i.e. l'unique quotient irréductible, de $X_{E,I}=\iota_{L_{E,I}}^H(\tau_{E,I})$.
\end{itemize}

Ainsi $\tau_{E,I}$ est une représentation irréductible essentiellement de carré intégrable de $L_{E,I}$. Observons que si $I\subset J$ alors $X_{E,J}$ est une sous-représentation de $X_{E,I}$ si $I\subset J$. 
De plus, $\pi_{E,J}$ est un sous-quotient de $X_{E,I}$ si et seulement si $I \subset J$. Les représentations $\pi_{E,I}$, qui sont donc les sous-quotients irréductibles de $$X_{E, \emptyset}=\nu^{k-1}\rho_E\times \nu^{k-2}\rho_E\times \dots \times \rho_E,$$ apparaissent avec multiplicité 1 dans la représentation $X_{E,\emptyset}$.

Alors, les représentations elliptiques de $H$ sont exactement les représentations $\pi_{E,I}$ ainsi construites à partir d'une représentation essentiellement de carré intégrable $\tau_E$ de $H$. 

Notons qu'une représentation irréductible de $H$ est elliptique si et seulement si elle a même support cuspidal qu'une représentation irréductible essentiellement de carré intégrable, en l'occurrence un segment de Zelevinski.

\vskip4mm
Nous avons la même construction pour les représentations elliptiques de $G=\text{GL}_n(F)$. 

\section{Résultats connus d'induction automorphe que l'on va utiliser}
Nous exposons ici les résultats déjà démontrés d'induction automorphe. Cela concerne les représentations essentiellement de carré intégrable.

\vskip3mm
Nous avons la proposition suivante dans \cite[p.148]{HLSMF}.

\prop \label{inductionautomorpheconnu} 
\begin{enumerate}
\item Soit $\tau_E$ une représentation irréductible cuspidale de $H$. Si la classe d'isomorphisme de $\tau_E$ a un stabilisateur d'ordre $d_1$ dans $\Gamma=\text{Gal}(E/F)$, alors son $\kappa$-relèvement $\pi$ est induite parabolique de $\pi_1\otimes \kappa\pi_1\otimes \dots \otimes \kappa^{d_1-1}\pi_1$ à $G$, où $\pi_1$ est une représentation irréductible cuspidale de $\text{GL}_{n_1}(F)$, $n=n_1d_1$, et a pour stabilisateur $\kappa^{d_1\mathbb{Z}}$ dans $\kappa^{\mathbb{Z}}$.

\item Si $\tau_E$ est essentiellement de carré intégrable, elle est déterminée par son support cuspidal qui forme un "segment" $\{\rho_E,\nu_E\rho_E,\ldots , \nu_E^{k-1}\rho_E\}$ (cf. 2.1), où $\rho_E$ est une représentation irréductible cuspidale de $\text{GL}_s(E), sk=m$, et $\nu_E=\nu\circ N_{E/F}$. D'après le point précédent on peut écrire le $\kappa$-relèvement de $\rho_E$ comme induite parabolique de $\pi_1\otimes \kappa\pi_1\otimes \dots \otimes \kappa^{d_1-1}\pi_1$ à $\text{GL}_{sd}(F), sd=n_1d_1$. Alors le $\kappa$-relèvement de $\tau_E$ est induite parabolique de $\pi_1'\otimes \kappa\pi_1'\otimes \dots \otimes \kappa^{d_1-1}\pi_1'$ à $G$, où $\pi_1'$ est la représentation essentiellement de carré intégrable de $\text{GL}_{n_1k}(F)$ de support cuspidal $\{\pi_1, \nu\pi_1, \dots ,\nu^{k-1}\pi_1\}$.

\end{enumerate}

\section{Compatibilité induction automorphe - induction parabolique}

Comme nous pouvons le voir dans la construction des représentations elliptiques l'induction parabolique est très présente. Nous nous intéressons donc à la question de la compatibilité entre l'induction parabolique et l'induction automorphe. 

\vskip3mm

Nous avons la proposition suivante dans \cite[p.145]{HLSMF}, on en donne les notations introduites. On se donne des entiers strictement positifs $m_1, \dots, m_t$ tels que $\sum_{i=1}^tm_i=m$. Pour $i=1, \dots, t$, on choisit un élément $e_i$ de $E^{\times}$ tel que $\sigma(e_i)=(-1)^{m_i(d-1)}e_i$, ce qui permet de considérer les facteurs de transfert $\tilde{\Delta}_i$ et $\Delta_i$ relatifs à l'induction automorphe de $H_i=\text{GL}_{m_i}(E)$ à $G_i=\text{GL}_{m_id}(F)$. Pour $i=1, \dots , t$ on se donne une base du $F$-espace vectoriel $E^{m_i}$, ce qui donne un plongement de $H_i$ dans $G_i$. Voyant $E^m$ comme $E^{m_1}\oplus \dots \oplus E^{m_t}$, on obtient une base du $F$-espace vectoriel $E^m$ d'où un plongement de $H$ dans $G$. Le groupe $L=G_1\times \dots \times G_t$ apparaît comme un sous-groupe de Levi de $G$, $L_H=H_1\times \dots \times H_t$ comme un sous-groupe de Levi de $H$, et on a $L_H=L\cap H$.

Soit $P$ le sous-groupe parabolique de $G$ formé des matrices triangulaires inférieures par blocs de taille $m_1d, \dots , m_td$, et soit $U_P$ le radical unipotent de $P$.

Le groupe $P_H=P \cap H$ est un sous-groupe parabolique de $H$, de radical unipotent $U_{P,H}=U_P\cap H$, et $L_H$ est une composante de Levi de $P_H$.

Pour $i=1, \dots, t$ on se donne une représentation $\pi_i$ de $G_i$.

\vskip3mm

\prop \label{compatibilitéparaboliqueautomorphe} 
Supposons que pour $i=1, \dots, t$, la représentation (irréductible, $\kappa$-stable) $\pi_i$ de $G_i$  soit un $\kappa$-relèvement d'une représentation lisse irréductible $\tau_i$ de $H_i$, et que les représentations $\pi=\iota_P^G(\pi_1\otimes \dots \otimes \pi_t)$ de $G$ et $\tau=\iota_{P_H}^H(\tau_1\otimes \dots \otimes \tau_t)$ de $H$ soient irréductibles. Alors $\pi$ est un $\kappa$-relèvement de $\tau$.
De plus, il existe une racine de l'unité $\zeta$, qui ne dépend ni des $\pi_i$, ni des $\tau_i$, telle que si pour $i=1, \dots, t$, $A_i$ est un isomorphisme de $\kappa\pi_i$ sur $\pi_i$, et que $A$ est l'isomorphisme de $\kappa\pi$ sur $\pi$ associé aux $A_i$ comme plus haut, on ait 
\[c(\tau, \pi, A)=\zeta \Pi_{i=1}^tc(\tau_i,\pi_i,A_i).\]

\section{Induction automorphe pour les représentations elliptiques}

On reprend dans cette section 7 les notations introduites dans la section 4. 
Le théorème suivant est le  résultat principal de l'article.

\theo Toute représentation irréductible elliptique de $H$ admet un $\kappa$-relèvement.

\vskip3mm

\noindent \textbf{Démonstration}

\begin{enumerate}

\item Nous partons donc d'une représentation essentiellement de carré intégrable $\tau_E$ de $H$ de support cuspidal $\{\rho_E, \nu_E\rho_E,\ldots , \nu_E^{k-1}\rho_E\}$ avec $k\vert m$. Comme nous l'avons vu dans la section 4, cette représentation $\tau_E$ permet de construire des représentations elliptiques $\pi_{E,I}$ de $H$ où $I$ est un sous-ensemble de $\mathcal{K}= \{1,\ldots ,k-1\}$.

Nous allons montrer que $\pi_{E,I}$ admet un  $\kappa$-relèvement.

D'après la proposition \ref{inductionautomorpheconnu}, nous savons qu'il existe des entiers $n_1, d_1$ avec $kn_1d_1=n$, que $\rho_E$ a un $\kappa$-relèvement de la forme (induite parabolique irréductible) 
$\rho \times \kappa\rho \times \cdots \times \kappa^{d_1-1}\rho$ pour une représentation irréductible cuspidale $\rho$ de ${\rm GL}_{n_1}(F)$ et qu'il existe une représentation $\xi$ de $\text{GL}_{n_1k}(F)$ (la représentation essentiellement de carré intégrable de ${\rm GL}_{kn_1}(F)$ de support cuspidal $\{\rho, \nu \rho , \ldots , \nu^{k-1}\rho\}$) tels que le $\kappa$-relèvement de $\tau_E$ soit de la forme \[\pi:= \xi\times \kappa\xi\times \dots \times \kappa^{d_1-1}\xi.\]

Montrons que le $\kappa$-relèvement de $\pi_{E,I}$ est \[\pi_I:=\sigma_I\times \kappa\sigma_I\times \dots \times \kappa^{d_1-1}\sigma_I\] où $\sigma_I$ est la représentation irréductible elliptique 
de ${\rm GL}_{kn_1}(F)$ associée à $\xi$ et $I$.

\vskip3mm

Pour cela nous allons montrer qu'il existe une constante $c$ telle que  pour toutes fonctions $f\in \mathcal{C}_c^{\infty}(H)$ et $\phi\in \mathcal{C}_c^{\infty}(G)$ qui se correspondent, on ait la relation \[\text{tr}\left(\pi_I(\phi)A_{\pi_I}\right)=c \text{ } \text{tr}\left(\pi_{E,I}(f)\right).\]

\vskip3mm
\item Nous introduisons dans ce paragraphe une représentation $\Theta$ telle que les $\pi_I$ définies ci-dessus en soient les sous-quotients irréductibles $\kappa$-stables. Cela permettra de déterminer plus facilement les opérateurs de $\kappa$-entrelacement associés aux $\pi_I$, opérateurs nécessaires pour montrer ce que l'on veut.

\vskip2mm

Soit $\rho$ la représentation cuspidale associée à $\xi$ via la classification de Bernstein-Zelevinsky.

Soit $\Theta$ la représentation induite 
\begin{equation*}
	\begin{split}
	(\nu^{k-1}\rho\times \nu^{k-1}\kappa\rho\times \dots \times \nu^{k-1}\kappa^{d_1-1}\rho ) \times (\nu^{k-2}\rho\times \nu^{k-2}\kappa\rho\times \dots \times \nu^{k-2}\kappa^{d_1-1}\rho )  \times  \dots \\
	\dots  \times (\rho\times\kappa\rho\times \dots \times \kappa^{d_1-1}\rho).
	\end{split}	
\end{equation*}

Cette représentation est isomorphe à 
\[\left(\nu^{k-1}\rho\times \nu^{k-2}\rho\times \dots \rho\right) \times \left(\nu^{k-1}\kappa\rho\times\dots \kappa\rho\right)\times \dots\times \left(\nu^{k-1}\kappa^{d_1-1}\rho\times \nu^{k-2}\kappa^{d_1-1}\rho\times \dots \times \kappa^{d_1-1}\rho\right). \]

En notant $\Theta_i=\nu^{k-1}\kappa^{i-1}\rho\times\dots \times \kappa^{i-1}\rho$ nous obtenons $\Theta=\Theta_1\times \Theta_2\times \dots \times \Theta_{d_1}$ et de plus grâce à la construction de la section 4 nous savons que $\Theta_i=\iota_{L_1}^{G_1}(\kappa^{i-1}\xi_{\emptyset})=\kappa^{i-1} \Theta_1$ avec $G_1= {\rm GL}_{kn_1}(F)$ et $L_1= {\rm GL}_{n_1}(F) \times \cdots \times {\rm GL}_{n_1}(F)$.

Or, nous connaissons les sous-quotients irréductibles de $\Theta_1$, ce sont précisément les $\sigma_I$ pour $I\subset \mathcal{K}$.

\vskip2mm
D'après \cite[Prop. 8.5]{Zelevinsky}, la représentation $\nu^a\kappa^{i}\rho\times \nu^b\kappa^{j}\rho$ est irréductible et isomorphe à $\nu^b\kappa^{j}\rho\times \nu^a\kappa^{i}\rho$ pour tous $0\leqslant i<j\leqslant d_1-1$ et tous $a,b$ entiers. Donc aucun sous-quotient irréductible de $\Theta_i$ n'est isomorphe à un sous-quotient irréductible de $\Theta_j$.

Donc, les sous-quotients irréductibles de $\Theta$ sont de multiplicité 1 et de la forme $\sigma_{I_1}\times \kappa\sigma_{I_2}\times \dots \times \kappa^{d_1-1}\sigma_{I_{d_1}}$, où $I_1, \dots, I_{d_1} \in \mathcal{P}(\mathcal{K})$.

\vskip2mm

Alors, pour $I\subset \mathcal{K}$, les $\pi_I= \sigma_I \times \kappa \sigma_I \times \cdots \times \kappa^{d_1-1}\sigma_I$ sont les sous-quotients irréductibles de $\Theta$ qui sont $\kappa$-stables.

\item  Déterminons maintenant les opérateurs de $\kappa$-entrelacement normalisés $A_{\pi_I}$ pour $I\subset \mathcal{K}$.

Pour chaque $I\subset \mathcal{K}$, notons  $\Xi_I= \iota_{L_{1,I}}^{G_1}(\xi_I)$  où $\xi_I$ est l'unique sous-représentation irréductible de l'induite parabolique $\iota_{L_1}^{L_{1,I}}(\nu^{k-1}\rho \otimes \cdots \nu \rho \otimes \rho)$. Nous avons donc $\sigma_I=L(\Xi_I)$.

Posons  \[\Xi_{(I)}=\Xi_I\times \kappa\Xi_I\times \dots \times \kappa^{d_1-1}\Xi_I.\]
Alors $\Xi_{(I)}$ est une sous-représentation de $\Theta$ et les sous-quotients irréductibles de $\Xi_{(I)}$ sont les $\sigma_{I_1}\times \kappa\sigma_{I_2}\times \dots \times \kappa^{d_1-1}\sigma_{I_{d_1}}$ avec $I\subset I_i$ pour chaque $i\in \lbrace 1, \dots , d_1\}$.

Nous remarquons que $\Xi_{(I)}$ est $\kappa$-stable. Donc $\Xi_{(I)}$ est $A_{\Theta}$-stable, où $A_{\Theta}$ est l'opérateur de $\kappa$-entrelacement obtenu grâce à l'induction parabolique à partir de $A_{\nu^{k-1}u}\otimes \dots \otimes A_u$ avec (rappel) $A_{\nu^i u}= A^{\rm g\acute{e}n}(\nu^i u, \psi)$.

Notons $\xi_I=\xi_I^1\otimes \dots \otimes\xi_I^{m_I}$ où $m_I$ est le nombre de blocs de $L_{1,I}$ et où les $\xi_I^1, \dots, \xi_I^{m_I}$ sont des représentations essentiellement de carré intégrable.

Notons, pour $1\leqslant j\leqslant m_I$, $\xi_{(I)}^j=\xi_I^j\times \kappa \xi_I^j\times \dots \times \kappa^{d_1-1}\xi_I^j$ (induite parabolique de $L_I= L_{1,I} \times \cdots \times L_{1,I}$ à $G$).

D'après \cite[prop 2.2, 2.3]{Tadic} nous avons \[\pi_I=L\left(\xi_{(I)}^1, \dots, \xi_{(I)}^{m_I}\right).\]

Pour $j=1,\ldots ,m_I$, notons $\alpha=$ $\alpha_I^j$ la longueur du segment de $\xi_I^j$. Alors, $\nu^{k-1}u\times \nu^{k-2}u\times \dots \times \nu^{k-\alpha}u$ est une sous-représentation d'une représentation induite à partir d'un segment de longueur $\alpha d_1$ et admet $\xi_{(I)}^{j}$ comme sous-quotient irréductible de multiplicité 1. Comme $\xi_{(I)}^{j}$ est générique, on peut lui appliquer la proposition \ref{compatibilitéinductionnormalisation} qui nous dit que son opérateur de $\kappa$-entrelacement normalisé est obtenu à partir de $A_{\nu^{k-1}u}\otimes \dots \otimes A_{\nu^{k-\alpha}u}$ par la propriété d'induction parabolique de multiplicité 1.

Nous concluons grâce à la proposition \ref{transitivité} que pour tout $I\subset \mathcal{K}$, $A_{\Theta}(\pi_I)=A_{\pi_I}$.

\item Fixons deux fonctions $\phi$ et $f$ qui se correspondent. Pour rappel, nous avons noté 
$\Xi_I$ la représentation $\iota_{L_{1,I}}^{G_1}(\xi_I)$ de $G_1= {\rm GL}_{\frac{n}{d_1}}(F)$, $\sigma_I$ le quotient de Langlands de $\Xi_I$ et $\Xi_{(I)}$ la représentation \[\Xi_{(I)}:=\Xi_I\times \kappa \Xi_I\times \dots \times \kappa^{d_1-1}\Xi_I.\]

\vskip3mm
Montrons maintenant que \[\text{tr}\left(\Xi_{(I)}(\phi)A_{\Xi_{(I)}}\right)=\sum_{I\subset J}\text{tr}\left(\pi_J(\phi)A_{\pi_J}\right).\]

Soit $0\subset U_1=\pi_{\mathcal{K}}\subset U_2\subset \dots \subset U_m=\Xi_{(I)}$ une suite de Jordan-Hölder pour l'action de $\text{GL}_n(F)$ via $\Theta$ et $A_{\Theta}$, i.e. tous les sous-modules dans la suite sont stables à la fois par $\Theta$ et $A_{\Theta}$ et que les quotients $U_{i+1}/U_i$ sont irréductibles pour cette action.

D'une part nous avons: \[\text{tr}\left( \Xi_{(I)}(\phi)A_{\Xi_{(I)}} \right)=\sum_{i=1}^m\text{tr}\left(U_{i+1}/U_i(\phi)A_{\Theta}(U_{i+1}/U_i)\right)\]

Or, si $U_{i+1}/U_i$ n'est pas irréductible alors $\text{tr}\left(U_{i+1}/U_i(\phi)A_{\Theta}(U_{i+1}/U_i)\right)=0$. 

En effet, si $U_{l+1}/U_l$ est un tel quotient, soit $\epsilon$ une sous-représentation irréductible pour l'action de $\text{GL}_n(F)$. Alors $\epsilon$ est isomorphe à une représentation de la forme $\sigma_{I_1}\times \kappa\sigma_{I_2}\times \dots \times \kappa^{d_1-1}\sigma_{I_{d_1}}$ avec les $I_i$ non tous égaux. Alors $A_{\Theta}$ envoie $\epsilon$ sur une autre sous-représentation irréductible de $\Xi_{(I)}$. 

Le quotient $U_{l+1}/U_l$ est la somme des conjugués de $\epsilon$ sous $A_{\Theta}$. Donc, s'il y a plus d'un conjugué et s'ils sont permutés par $A_{\Theta}$ sans point fixe, alors la trace est nulle.

\vskip2mm
Il ne reste donc dans la trace que les représentations irréductibles, à savoir les $\pi_J$ pour $I\subset J$: 

$$\text{tr}\left( \Xi_{(I)}(\phi)A_{\Xi_{(I)}} \right)=\sum_{I\subset J}\text{tr}\left(\pi_J(\phi)A_{\Theta}(\pi_J)\right).$$

Or, d'après le paragraphe précédent, $A_{\Theta}(\pi_J)=A_{\pi_J}$. D'où 
\[\text{tr}\left(\Xi_{(I)}(\phi)A_{\Xi_{(I)}}\right)=\sum_{I\subset J}\text{tr}\left(\pi_J(\phi)A_{\pi_J}\right).\]

\vskip3mm

\item  De plus, par compatibilité de l'application de $\kappa$-relèvement avec l'induction parabolique (proposition \ref{compatibilitéparaboliqueautomorphe}), $\Xi_{(I)}=\Xi_I\times \kappa\Xi_I \times \dots \times \kappa^{d_1-1}\Xi_I$ est 
un $\kappa$-relèvement de $X_{E,I}$ où $X_{E,I} = \iota_{L_{E,I}}^{H}(\tau_{E,I})$.

Il existe donc une constante $c\in \mathbb{C}$ telle que $\text{tr}\left(\Xi_{(I)}(\phi)A_{\Xi_{(I)}}\right)=c ~\text{tr}\left(X_{E,I}(f)\right)$.

\item Or, sur $\text{GL}_m(E)$ nous avons \[\text{tr}X_{E,I}(f)=\sum_{I\subset J}\text{tr}\pi_{E,J}(f).\]

\item Nous avons donc $$c ~\text{tr}\left(X_{E,I}(f)\right)= \text{tr}\left(\Xi_{(I)}(\phi)A_{\Xi_{(I)}}\right)= \sum_{I\subset J}\text{tr}\left(\pi_J(\phi)A_{\pi_J}\right)$$ et $$\text{tr}X_{E,I}(f)=\sum_{I\subset J}\text{tr}\pi_{E,J}(f).$$

Donc par récurrence décroissante nous obtenons \[\text{tr}\left(\pi_I(\phi)A_{\pi_I}\right)=c ~ \text{tr}\left(\pi_{E,I}(f)\right).\]

\end{enumerate}

Cela achève la démonstration du théorème. \hfill \qed

\vskip2mm

\printbibliography[]

\end{document}